# NUMERICAL SOLUTION OF CONSERVATIVE FINITE-DIMENSIONAL STOCHASTIC SCHRODINGER EQUATIONS[1]

By Carlos M. Mora

*Universidad de Concepción*

The paper deals with the numerical solution of the nonlinear Itô stochastic differential equations (SDEs) appearing in the unravelling of quantum master equations. We first develop an exponential scheme of weak order 1 for general globally Lipschitz SDEs governed by Brownian motions. Then, we proceed to study the numerical integration of a class of locally Lipschitz SDEs. More precisely, we adapt the exponential scheme obtained in the first part of the work to the characteristics of certain finite-dimensional nonlinear stochastic Schrödinger equations. This yields a numerical method for the simulation of the mean value of quantum observables. We address the rate of convergence arising in this computation. Finally, an experiment with a representative quantum master equation illustrates the good performance of the new scheme.

## 1. Introduction.

1.1. *Objectives.* The primary objective of this paper is to develop an efficient scheme for the computation of $\mathbf{E}\langle Z_t, AZ_t\rangle$, where $A \in \mathbb{C}^{d,d}$, $\langle \cdot, \cdot \rangle$ is the standard scalar product in $\mathbb{C}^d$, and $Z_t$ satisfies the following Itô stochastic differential equation (SDE) on $\mathbb{C}^d$:

$$(1) \qquad Z_t = Z_0 + \int_0^t (GZ_s + D(Z_s))\,ds + \sum_{k=1}^n \int_0^t E_k(Z_s)\,dW_s^k,$$

Received March 2004; revised December 2004.
[1]Supported in part by FONDECYT Grants 1030552 and 1040623, Grant Milenio ICM P02-049 and the International scientific cooperation program CONICYT/INRIA.

where $\|Z_0\| = 1$, $W$ is an $n$-dimensional $(\mathfrak{F}_t)$-Brownian motion,

$$G = -iH - \tfrac{1}{2} \sum_{k=1}^{n} L_k^* L_k,$$

(2)

$$D(z) = \sum_{k=1}^{n} (\operatorname{Re}\langle z, L_k z\rangle L_k z - \tfrac{1}{2} \operatorname{Re}^2 \langle z, L_k x\rangle z),$$

provided $z \in \mathbb{C}^d$, and for any $k = 1, \ldots, n$ and $z \in \mathbb{C}^d$

(3) $$E_k(z) = L_k z - \operatorname{Re}\langle z, L_k z\rangle z.$$

Here $L_k \in \mathbb{C}^{d,d}$ for all $k = 1, \ldots, n$, $H$ is a $d \times d$ self-adjoint matrix, and $(\Omega, \mathfrak{F}, \mathbf{P}, (\mathfrak{F}_t)_{t \geq 0})$ is the underlaying filtered probability space. Our main motivation came from the numerical simulation of the evolution of open quantum systems. To help to shed light on our problem, Section 1.2 below looks closely at this application.

In this work we follow the strategy of constructing exponential schemes adapted to the characteristics of (1). This approach is partially motivated by the good behavior of the exponential integrators in the solution of certain class of stiff ordinary differential equations, for example, those associated to both time-dependent Schrödinger equations and oscillatory electric circuits (see, e.g., [12, 13] for more details). Another motivation came from a number of numerical experiments which illustrate the good performance of the exponential schemes for real SDEs with additive noise whose numerical solution by the standard explicit schemes presents numerical instabilities (see, e.g., [4, 16, 24]).

To gain understanding of the exponential methods, this article starts by generalizing the Euler-exponential scheme for SDEs with additive noise proposed in [24] to the context of SDEs of the form

(4) $$X_t = X_0 + \int_0^t b(s, X_s)\, ds + \int_0^t \sigma(s, X_s)\, dW_s,$$

where $t \in [0, T]$, $X_t$ takes values in $\mathbb{R}^d$ and $b$, $\sigma$ are smooth functions with bounded derivatives up to appropriate order. Indeed, adapting the methodology employed in [24], we develop the following numerical method:

SCHEME 1 (Euler-exponential).  Let $\xi_0^1, \ldots, \xi_0^n, \xi_1^1, \ldots, \xi_1^n, \ldots, \xi_{M-1}^1, \ldots, \xi_{M-1}^n$ be independent and identically distributed (i.i.d.) real random variables with symmetric law, variance 1 and moments of any order. Then, we consider the recursive algorithm

$$V_{m+1} = \exp\!\left(Jb(T_m, V_m)\frac{T}{M}\right)$$
$$\times \left(V_m + \frac{T}{M}(b(T_m, V_m) - Jb(T_m, V_m)V_m) + \sqrt{\frac{T}{M}}\sigma(T_m, V_m)\xi_m\right),$$

where $T_m = mT/M$, $Jb = (\partial_x^j b^k)_{k,j=1,\ldots,d}$ and $\xi_m = (\xi_m^1, \ldots, \xi_m^n)^\top$. Here $V_0$ is independent of $\xi_m$ for all $m = 1, \ldots, M-1$.

In particular, we check that the error between $\mathbf{E} f(X_T)$ and $\mathbf{E} f(V_M^M)$ has a linear behavior when $f$ is a smooth function, that is, Scheme 1 achieves the first order of weak convergence. Preliminary numerical experiments suggest that Scheme 1 should be useful in situations where the eigenvalues of $Jb$ have vastly different sizes and their real parts are nonpositive. In these circumstances, both the explicit Euler scheme and the implicit Euler scheme present relevant time step restrictions in many cases. It will be interesting to test more carefully Scheme 1 with theoretical and real-life scientific problems.

Combining Scheme 1 with splitting and projection techniques yields an efficient numerical method for (1). To be more precise, we may split the drift term of (1) into two components to obtain that there exists a continuous semimartingale $S_{\cdot,T_m}$ such that

$$\text{(5)} \qquad Z_t = Z_{T_m} + \int_{T_m}^t G Z_r \, dr + S_{t,T_m}$$

for all $t \in [T_m, T_{m+1}]$. Here $T_m = mT/M$. Then, solving explicitly, the linear SDE (5) leads to

$$\text{(6)} \qquad Z_t = \exp(G(t - T_m)) Z_{T_m} + \int_{T_m}^t \exp(G(t - r)) \, dS_{r,T_m},$$

for any $t \in [T_m, T_{m+1}]$. Letting $\hat{Z}_m \approx Z_{T_m}$ and replacing the right-hand side of (6) by random vectors with similar first three moment properties, we arrive at the weak approximation $\Phi_{m+1}^{\hat{Z}_m, M}$ of $Z_{T_{m+1}}$, where

$$\text{(7)} \qquad \Phi_{m+1}^{z, M} = \exp\left(G \frac{T}{M}\right) \left( z + D(z) \frac{T}{M} + \sqrt{\frac{T}{M}} \sum_{k=1}^n E_k(z) \xi_m^k \right),$$

with $\xi_0^1, \ldots, \xi_0^n, \ldots, \xi_{M-1}^1, \ldots, \xi_{M-1}^n$ as in Scheme 1. Now, the growth of $\Phi_{m+1}^{\hat{Z}_m, M}$ is stabilized by using a projection technique from the numerical treatment of ordinary differential equations with invariants. To be precise, since $\|Z_t\| = 1$, we project $\Phi_{m+1}^{\hat{Z}_m, M}$ onto the surface of the unit ball. From this we derive the following numerical method:

SCHEME 2 (Version of the Euler-exponential method). Let $\hat{Z}_0^M$ be a random variable independent of $\xi_1, \ldots, \xi_{M-1}$ satisfying $\|\hat{Z}_0^M\| = 1$. Then, we define recursively

$$\hat{Z}_{m+1}^M = p\left( \exp\left(G \frac{T}{M}\right) \left( \hat{Z}_m^M + D(\hat{Z}_m^M) \frac{T}{M} + \sqrt{\frac{T}{M}} \sum_{k=1}^n E_k(\hat{Z}_m^M) \xi_m^k \right) \right),$$

where

$$p(z) = \begin{cases} 0, & \text{if } z = 0, \\ z/\|z\|, & \text{if } z \neq 0. \end{cases}$$

The behavior of $\hat{Z}_m$ is tested by means of a numerical experiment. Sometimes the application of the projection techniques deteriorates the performance of the numerical method. For instance, Section VII.2 of [10] presents an example where a projection procedure destroys the correct qualitative behavior of a symplectic Euler method applied to a deterministic Hamiltonian system. In our case, a numerical experiment with a representative quantum system, where the eigenvalues of $G$ have vastly different sizes and their real parts are nonpositive, illustrates the very good behavior of $\hat{Z}_m$. In this example, both versions explicit and implicit of the Euler scheme fail.

A secondary objective of this paper is to study convergence properties of the stochastic schemes used for computing $\mathbf{E}\langle Z_t, AZ_t\rangle$. In particular, our interest is focused on the rate of convergence of $\mathbf{E}\langle \hat{Z}_m, A\hat{Z}_m\rangle$. Most of the existing convergence theory for numerical methods requires that the coefficients of the SDE be globally Lipschitz. This motivates the increasing interest in addressing convergence properties of the numerical schemes for more general class of SDEs (see, e.g., [11, 20, 31]). In our case, the SDE under, consideration is only locally Lipschitz. In fact, the coefficients of (1) have nonlinear grow. To deal with this situation, we modify the standard arguments due to Talay [27, 28] and Milshtein [21]. Another difficulty in carrying out our theoretical study is that $p$ has a singularity at 0. To overcome it, we take some inspirations in [29] and [9].

This paper is organized in six sections. Section 2 is devoted to introduce notation. Section 3 develops Scheme 1. In Section 4 we construct Scheme 2. Section 5 provides the rate of convergence of $\mathbf{E}\langle \hat{Z}_M^M, A\hat{Z}_M^M\rangle$. Section 6 presents a numerical experiment.

1.2. *Motivation.* We start with a brief exposition of some basic results of quantum mechanics. The states of a quantum system are described by elements of an adequate complex Hilbert space $(\mathfrak{h}, \langle \cdot, \cdot \rangle)$ and the quantum observables are represented by self-adjoint linear operators in $\mathfrak{h}$ (see, e.g., [5] for more details). In the Heisenberg picture, the evolution of the observable $A$ under the Born–Markov approximation is given by the minimal solution of the adjoint quantum master equation

$$(8) \qquad \frac{d}{dt}\tau_t = G^*\tau_t + \tau_t G + \sum_{k=1}^{n} L_k^* \tau_t L_k, \qquad \tau_0 = A.$$

Here $\tau_t, L_1, \ldots, L_n$ are general linear operators in $\mathfrak{h}$ and

$$G = -iH - \tfrac{1}{2}\sum_{k=1}^{n} L_k^* L_k,$$

with $H$ self-adjoint operator in $\mathfrak{h}$. The operators $L_j$, with $j = 1, \ldots, n$, describe the effects of the environment and $H$ represents the Hamiltonian. For a fuller mathematical treatment of (8), we refer the reader to [8].

The concept of quantum trajectories allows the simulation of the evolution of the quantum observables. Let us consider the linear stochastic evolution equation on $\mathfrak{h}$,

$$(9) \qquad Y_t = y_0 + \int_0^t GY_s \, ds + \sum_{k=1}^n \int_0^t L_k Y_s \, dB_s^k,$$

where $B$ is an $n$-dimensional Brownian motion on the filtered complete probability space $(\Omega, \mathfrak{F}, \mathbf{Q}, (\mathfrak{F}_t)_{t \geq 0})$. Then, in the measurement interpretation of the quantum trajectories the stochastic process $Z_t = Y_t / \|Y_t\|$ describes the evolution of a system conditioned on continuous observation (see, e.g., [33] and the references given there). Furthermore, the mean value of the observable $A$ at the instant $t$ is given by $\mathbf{E}_\mathbf{Q} \langle Y_t, AY_t \rangle$. In fact, we may see that

$$(10) \qquad \mathbf{E}_\mathbf{Q} \langle Y_t, AY_t \rangle = \langle y_0, \tau_t y_0 \rangle,$$

under certain assumptions (see, e.g., [3, 15]).

The stochastic Schrödinger equations (1) allows the description of finite-dimensional open quantum systems, for example, $q$-bits models. Let $\dim \mathfrak{h} < +\infty$. Applying Itô's formula, integration by parts formula and Girsanov's theorem, we see that there exists a probability measure $\mathbf{P}$, which is equivalent to $\mathbf{Q}$, such that $Z_t$ satisfies (1) and

$$(11) \qquad \mathbf{E}_\mathbf{Q} \langle Y_t, AY_t \rangle = \|y_0\|^2 \mathbf{E}_\mathbf{P} \langle Z_t, AZ_t \rangle.$$

Combining (10) and (11) gives

$$\langle y_0, \tau_t y_0 \rangle = \|y_0\|^2 \mathbf{E}_\mathbf{P} \langle Z_t, AZ_t \rangle.$$

Hence, the numerical solution of (1) leads to the numerical description of $\langle y_0, \tau_t y_0 \rangle$, which represents the mean value of the observable $A$ at the instant $t$. This procedure has been proposed in the physical literature in order to overcome the difficulties appearing in the direct numerical integration of (8) (see, e.g., [25]). It is worth pointing out that the numerical schemes for (8) exhibit serious numerical instabilities and the dimension of the state space of (8) grows up very fast to $+\infty$ when $\dim \mathfrak{h} \to +\infty$. On the other hand, the computation $\langle y_0, \tau_t y_0 \rangle$ by means of the numerical solution of (9) presents drawbacks. In this case, a large number of numerical experiments show the blow-up of the trajectories of the explicit Euler method even for small size of the discretization step. Furthermore, in many situations the implicit Euler scheme, defined as in [17, 20, 22, 31], tends to the origin very fast.

Finally, the efficient numerical solution of (1) also plays an important role in the study of many infinite-dimensional quantum phenomena, for instance, harmonic oscillators. Let $\dim \mathfrak{h} = +\infty$. Then, we may choose an adequate sequence $(\mathfrak{h}_d)_d$ of finite-dimensional Hilbert spaces such that $\mathbf{E}_{\mathbf{Q}}\langle Y_t, AY_t \rangle$ is approximated by $\mathbf{E}_{\mathbf{Q}}\langle Y_{t,d}, AY_{t,d} \rangle$, where $Y_{t,d}$ is the continuous adapted stochastic process with values on $\mathfrak{h}_d$ given by

$$Y_{t,d} = P_d Y_0 + \int_0^t G_d Y_{s,d}\, ds + \sum_{k=1}^n \int_0^t P_d L_k Y_{s,d}\, dB_s^k,$$

with $P_d : \mathfrak{h} \to \mathfrak{h}_d$ the orthogonal projection of $\mathfrak{h}$ over $\mathfrak{h}_d$ and

$$G_d = -iP_d H - \tfrac{1}{2} \sum_{j=1}^n P_d L_j^* P_d L_j$$

(see, e.g., [23], where the rate of convergence of this approximation is studied). Thus, similar arguments to those used in the previous paragraph give rise to our main problem.

**2. Notation.** Throughout this paper we assume that the filtered probability spaces satisfy the usual hypotheses (see, e.g., [7, 26]). We will denote by $\mathbf{E}_t$ the conditional expectation with respect to $\mathfrak{F}_t$. For simplicity, we restrict our attention to the equidistant partitions of the time interval $[0, T]$, that is, time discretizations of the form $(T_m^M)_{m=0,\ldots,M}$, with $T_m^M = mT/M$. To shorten notation, sometimes the explicit dependence on the discretization step $T/M$ will be suppressed except where we wish to emphasize its role. For example, we will write $T_m$ instead of $T_m^M$ if no misunderstanding is possible. We will use the same symbol $K(\cdot)$ (resp., $K$ and $q$) for different positive increasing functions (resp., positive real numbers) having the common property to be independent of $M$. Furthermore, it is understood that $q$ is greater than or equal to 2.

Let $A \in \mathbb{C}^{l,q}$. Then, the symbol $A^\top$ will stand for the transpose of $A$. Furthermore, $A^{k,j}$ will be the $(k,j)$th component of the matrix $A$ and $\|A\| = \sqrt{\sum_{k=1}^l \sum_{j=1}^q |A^{k,j}|^2}$. For any $x, y \in \mathbb{C}^d$, we will write $\langle x, y \rangle = \sum_{k=1}^d \overline{x^k} y^k$ and $\bar{x} = (\overline{x^1}, \ldots, \overline{x^d})$.

For each $l \in \mathbb{N}$, we define $\mathcal{P}_l$ to be $\{1, \ldots, d\}^l$. For any $\vec{p} \in \mathcal{P}_l$, with $l \in \mathbb{N}$, and $x \in \mathbb{T}^d$, with $\mathbb{T} \in \{\mathbb{R}, \mathbb{C}\}$, we set

$$F_{\vec{p}}(x) = \prod_{k=1}^l x^{p^k}$$

and $\partial_x^{\vec{p}} g(x) = \frac{\partial}{\partial x^{p^1}} \cdots \frac{\partial}{\partial x^{p^l}} g(x)$, provided that $g : \mathbb{T}^d \to \mathbb{T}$ is smooth enough. The symbol $\partial^0$ stands for the identity operator and $\mathcal{P}_0 = \{0\}$. We say that a family of functions $(f_\theta : [0, T] \times \mathbb{R}^d \to \mathbb{R})_{\theta \in \Theta}$ belongs to $\mathcal{C}_p^L([0, T] \times \mathbb{R}^d, \mathbb{R})$ if for any $\vec{p} \in \mathcal{P}_l$, with $l \leq L$,

(i) $\partial_x^{\vec{p}} f_\theta$ is a continuous function whenever $\theta \in \Theta$,
(ii) $|\partial_x^{\vec{p}} f_\theta(t, x)| \leq K(T)(1 + \|x\|^q)$ for all $t \in [0, T]$, $x \in \mathbb{R}^d$, and $\theta \in \Theta$.

Moreover, we write $(f_\theta : \mathbb{R}^d \to \mathbb{R})_{\theta \in \Theta} \in \mathcal{C}_p^L(\mathbb{R}^d, \mathbb{R})$ if $(f_\theta)_{\theta \in \Theta}$ satisfies the conditions (i) and (ii).

**3. Euler-exponential scheme for general SDEs.** To shed some new light on exponential schemes, this section develops an exponential method for (4). First, we introduce a local exponential representation of $X$. Second, we derive a first weak order exponential scheme. Finally, we deal with the convergence analysis of the new scheme.

3.1. *Euler-exponential scheme.* Throughout this section, we assume that $b$ and $\sigma$ satisfy the standard conditions for the existence and uniqueness of strong solutions of (4) (see, e.g., [1, 26]). Furthermore, for any $s \in [0, T]$ and $x \in \mathbb{R}^d$, we consider the adapted stochastic process $X_t^{s,x}$ defined by

$$(12) \qquad X_t^{s,x} = x + \int_s^t b(r, X_r^{s,x}) \, dr + \int_s^t \sigma(r, X_r^{s,x}) \, dW_r,$$

for all $t \in [s, T]$. For abbreviation, we set $X_t^x := X_t^{0,x}$.

The following theorem states a local representation of $X$ of exponential type.

LEMMA 1. *Let $\partial_x^{\vec{p}} b$ be a continuous function for each $\vec{p} \in \mathcal{P}_3$. Suppose that $\partial_t b$ and $\partial_t \partial_x^k b$, with $k = 1, \ldots, d$, are also continuous functions. If $\xi$ is a $\mathfrak{F}_{T_m}$-random variable taking values in $\mathbb{R}^d$, then for any $t \in [T_m, T_{m+1}]$,*

$$\begin{aligned}
X_t^{T_m, \xi} &= e^{Jb(T_m, \xi)(t - T_m)} \xi + \int_{T_m}^t e^{Jb(T_m, \xi)(t-s)} \sigma(s, X_s^{T_m, \xi}) \, dW_s \\
(13) \qquad &+ \int_{T_m}^t e^{Jb(T_m, \xi)(t-s)} (b(T_m, \xi) - Jb(T_m, \xi)\xi) \, ds \\
&+ \int_{T_m}^t e^{Jb(T_m, \xi)(t-s)} \left( \int_{T_m}^s L(b)(u, X_u^{T_m, \xi}) \, du \right) ds + R_{t, T_m}.
\end{aligned}$$

*Recall from Section 1 that $Jb(T_m, \xi) = (\partial_x^j b^k(T_m, \xi))_{k,j=1,\ldots,d}$. In addition,*

$$L = \frac{\partial}{\partial t} + \frac{1}{2} \sum_{k,l=1}^d (\sigma \sigma^\top)^{k,l} \partial_x^{k,l}$$

*and $R_{t, T_m}$ is given by*

$$\sum_{j=1}^d \sum_{i=0}^n \int_{T_m}^t e^{Jb(T_m, \xi)(t-s)} \int_{T_m}^s \left( \int_{T_m}^r L_i(\partial_x^j b)(u, X_u^{T_m, \xi}) \, dW_u^i \right) b^j(r, X_r^{T_m, \xi}) \, dr \, ds$$

$$+ \sum_{j=1}^{d} \sum_{i=0}^{n} \int_{T_m}^{t} e^{Jb(T_m,\xi)(t-s)}$$

$$\times \int_{T_m}^{s} \left( \int_{T_m}^{r} L_i(\partial_x^j b)(u, X_u^{T_m,\xi}) \, dW_u^i \right) \sigma^{j,\cdot}(r, X_r^{T_m,\xi}) \, dW_r \, ds,$$

where

$$L_0 = \frac{\partial}{\partial t} + \sum_{k=1}^{d} b^k \, \partial_x^k + \frac{1}{2} \sum_{k,l=1}^{d} (\sigma \sigma^\top)^{k,l} \, \partial_x^{k,l}$$

and for any $i = 1, \ldots, n$,

$$L_i = \sum_{k=1}^{d} \sigma^{k,i} \, \partial_x^k.$$

PROOF. Using Itô's formula (see, e.g., [6]), we obtain, for any $s \geq T_m$,

$$b(s, X_s^{T_m,\xi}) = b(T_m, \xi) + \int_{T_m}^{s} L(b)(u, X_u^{T_m,\xi}) \, du$$

$$+ \sum_{j=1}^{d} \int_{T_m}^{s} \partial_x^j b(u, X_u^{T_m,\xi}) \, d(X_u^{T_m,\xi})^j.$$

Substituting this result into (12), we see that, for any $t \geq T_m$,

$$X_t^{T_m,\xi} = \xi + \int_{T_m}^{t} b(T_m, \xi) \, ds + \int_{T_m}^{t} \left( \int_{T_m}^{s} L(b)(u, X_u^{T_m,\xi}) \, du \right) ds$$

$$+ \sum_{j=1}^{d} \int_{T_m}^{t} \left( \int_{T_m}^{s} \partial_x^j b(u, X_u^{T_m,\xi}) \, d(X_u^{T_m,\xi})^j \right) ds + \int_{T_m}^{t} \sigma(s, X_s^{T_m,\xi}) \, dW_s.$$

Then, applying Itô's formula to each $\partial_x^i b$, with $i = 1, \ldots, d$, yields

$$(14) \quad X_t^{T_m,\xi} = \xi + \int_{T_m}^{t} Jb(T_m, \xi) X_s^{T_m,\xi} \, ds + S_{t,T_m} \qquad \forall t \in [T_m, T_{m+1}],$$

where

$$S_{t,T_m} = \int_{T_m}^{t} (b(T_m, \xi) - Jb(T_m, \xi)\xi) \, ds + \int_{r_m}^{t} \sigma(s, X_s^{T_m,\xi}) \, dW_s$$

$$+ \int_{T_m}^{t} \left( \int_{T_m}^{s} L(b)(u, X_u^{T_m,\xi}) \, du \right) ds$$

$$+ \sum_{j=1}^{d} \sum_{i=0}^{n} \int_{T_m}^{t} \left( \int_{T_m}^{r} \left( \int_{T_m}^{s} L_i(\partial_x^j b)(u, X_u^{T_m,\xi}) \, dW_u^i \right) b^j(s, X_s^{T_m,\xi}) \, ds \right) dr$$

$$+ \sum_{j=1}^{d} \sum_{i=0}^{n} \int_{T_m}^{t} \left( \int_{T_m}^{r} \left( \int_{T_m}^{s} L_i(\partial_x^j b)(u, X_u^{T_m,\xi}) \, dW_u^i \right) \sigma^{j,\cdot}(s, X_s^{T_m,\xi}) \, dW_s \right) dr,$$

with $W_u^0 := u$. Since $S_{\cdot,T_m}$ is a continuous semimartingale, the linear SDE (14) has the explicit solution established in (13) (see, e.g., [26]). □

We are now in position to deduce an exponential scheme of weak order 1 for (4). Since $X_{T_{m+1}} = X_{T_{m+1}}^{T_m, X_{T_m}}$, we start by replacing $\xi$ by $X_{T_m}$ in (13). Then, we neglect the terms of the last line of (13) because they involve multiple integrals. Applying the Itô–Taylor formula to $\sigma(s, X_s)$, we see that $\sigma(s, X_s)$ can be approximated by $\sigma(T_m, X_{T_m})$ in the first line of (13). Furthermore, we substitute $X_{T_m}$ for $\bar{X}_{T_m}$ in (13), where $\bar{X}_{T_m}$ is so chosen that it approximates $X_{T_m}$ in a weak linear sense. From this we may conclude that, for any $t \in [T_m, T_{m+1}]$,

$$X_t \approx Y_t^{\bar{X}_{T_m}, M},$$

where for each $\mathfrak{F}_{T_m}$-random variable $\xi$,

(15)
$$\begin{aligned}Y_t^{\xi,M} &= \exp(Jb(T_m,\xi)(t-T_m))\xi \\ &\quad + \int_{T_m}^t \exp(Jb(T_m,\xi)(t-s))(b(T_m,\xi) - Jb(T_m,\xi)\xi)\,ds \\ &\quad + \int_{T_m}^t \exp(Jb(T_m,\xi)(t-s))\sigma(T_m,\xi)\,dW_s.\end{aligned}$$

In order to replace $Y_{T_{m+1}}^{\bar{X}_{T_m}, M}$ by other random variables with similar first three moments, we use arguments similar to those in Section 3.2 of [24]. In fact, we approximate the integral

$$\int_{T_m}^{T_{m+1}} \exp(Jb(T_m,\bar{X}_{T_m})(t-s))(b(T_m,\bar{X}_{T_m}) - Jb(T_m,\bar{X}_{T_m})\bar{X}_{T_m})\,ds$$

by a classical rectangle. Finally, we look for a linear function $H_M(T_m, \cdot, \cdot)$ such that

$$\mathbf{E}_{T_m}(H_M(T_m, \bar{X}_{T_m}, \xi_m) H_M(T_m, \bar{X}_{T_m}, \xi_m)^\top),$$

with $\xi_0, \ldots, \xi_{M-1}$ defined as in Scheme 1, is the approximation of

$$\int_{T_m}^{T_{m+1}} e^{Jb(T_m,\bar{X}_{T_m})(T_{m+1}-s)} \sigma(T_m,\bar{X}_{T_m}) \sigma(T_m,\bar{X}_{T_m})^\top e^{Jb(T_m,\bar{X}_{T_m})^\top (T_{m+1}-s)}\,ds$$

given by a classical rectangle rule. This yields Scheme 1 defined in the Introduction, that is, the method

$$\begin{aligned}V_{m+1} &= \exp\left(Jb(T_m,V_m)\frac{T}{M}\right) \\ &\quad \times \left(V_m + \frac{T}{M}(b(T_m,V_m) - Jb(T_m,V_m)V_m) + \sqrt{\frac{T}{M}}\sigma(T_m,V_m)\xi_m\right).\end{aligned}$$

REMARK 1. Let the eigenvalues of the $Jb$ have nonpositive real part. Since
$$\exp\left(Jb(T_m, V_m)\frac{T}{M}\right) = \left(I - Jb(T_m, V_m)\frac{T}{M}\right)^{-1} + O\left(\left(Jb(T_m, V_m)\frac{T}{M}\right)^2\right),$$
Scheme 1 leads to the following version of the implicit Euler scheme which avoid the solution of nonlinear equations systems:
$$V_{m+1}^1 = \left(I - \frac{T}{M}Jb(T_m, V_m^1)\right)^{-1}$$
$$\times \left(V_m^1 + \frac{T}{M}(b(T_m, V_m^1) - Jb(T_m, V_m^1)V_m^1) + \sqrt{\frac{T}{M}}\sigma(T_m, V_m^1)\xi_m\right).$$

REMARK 2. To carry out the computation of $\exp(Jb(T_m, V_m)T/M)$ times a vector $u$ in the implementation of $V$, we may use Krylov approximations with Lanczos process (see, e.g., [12, 14]). In fact, many numerical experiments with high-dimensional problems illustrate the good performance of this numerical method. Furthermore, Hochbruck and Lubich [12] proved that the convergence of Krylov methods for $\exp(Jb(T_m, V_m)T/M)u$ is faster than that for the solution of the linear equation $(I - Jb(T_m, V_m)T/M)x = u$, which is required in both methods $V^1$ and the usual implicit Euler scheme when $b$ is linear. Alternative methods are Padé approximations, Strang splitting, Chebyshev approximations and Magnus integrators.

3.2. *Rate of convergence.* Similar to the Euler scheme (see, e.g., [2, 32]), the error between $\mathbf{E}f(X_T)$ and $\mathbf{E}f(V_M^M)$ can be expanded in powers of $T/M$ under general enough conditions. In particular, adapting the arguments used in [24] for studying the rate of weak convergence of the Euler-exponential scheme for SDEs with additive noise, we obtain the next theorem.

THEOREM 1. *Let $\mathbf{E}\|X_0\|^q < +\infty$ for any $q \geq 2$. Assume that $b$ and $\sigma$ are continuous functions such that $\partial_x^{\vec{p}}b$ and $\partial_x^{\vec{p}}\sigma$ are bounded continuous functions for every $\vec{p} \in \mathcal{P}_l$, with $l = 1, \ldots, 9$. Furthermore, suppose that the components of $\frac{\partial}{\partial t}\partial_x^{\vec{p}}b$, $\partial_x^{\vec{p}}\frac{\partial}{\partial t}b$, $\frac{\partial}{\partial t}\partial_x^{\vec{p}}\sigma$ and $\partial_x^{\vec{p}}\frac{\partial}{\partial t}\sigma$ belong to $\mathcal{C}_p^0([0,T] \times \mathbb{R}^d, \mathbb{R})$ whenever $\vec{p} \in \mathcal{P}_l$, with $l = 0, 1, 2$. Let the components of $\frac{\partial^2}{\partial t^2}b$, $\frac{\partial^2}{\partial t^2}\sigma$ belong to $\mathcal{C}_p^0([0,T] \times \mathbb{R}^d, \mathbb{R})$. If $f \in \mathcal{C}_p^9(\mathbb{R}^d, \mathbb{R})$, then there exists a continuous function $\Psi$, with $(\Psi(s, \cdot))_{s \in [0,T]} \in \mathcal{C}_p^4(\mathbb{R}^d, \mathbb{R})$, such that*

$$(16) \quad \left|\mathbf{E}f(X_T) - \mathbf{E}f(V_M^M) - \frac{T}{M}\int_0^T \mathbf{E}\Psi(s, X_s)\,ds\right|$$
$$\leq K(T)(1 + \mathbf{E}\|X_0\|^q)\left(\frac{T}{M}\right)^2,$$

*provided that $V_0$ have the same distribution as $X_0$.*

Before we prove Theorem 1, we present a series of observations and results. We start by considering the probability space $(\bar{\Omega}, \mathfrak{G}, \bar{\mathbf{P}})$ that arises from the completion of the product measure space induced by $(\Omega, \mathfrak{F}, \mathbf{P})$, $V_0$, and the random variables $\xi_0, \ldots, \xi_{M-1}$. Then, we combine the completion of $(\mathfrak{F}_t \otimes \sigma(V_0^M, \xi_k : k \leq [tM/T] - 1))_{t \geq 0}$ with a limit procedure to construct the filtration $(\mathfrak{G}_t)_{t \geq 0}$ that satisfies the usual hypotheses (see, e.g., IV 48 of [7]). Let us use from now on the same letter to designate a random variable and its natural extension to the Cartesian product space $\bar{\Omega}$, for instance, $(W_t)_{t \geq 0}$ also denotes the stochastic process $(W_t \circ \mathrm{Pr}_\Omega)_{t \geq 0}$, where $\mathrm{Pr}_\Omega$ is the projection of $\bar{\Omega}$ over $\Omega$. Thus, $W$ is an $n$-dimensional $(\mathfrak{G}_t)$-Brownian motion, $X_0$ and $V_0$ are $\mathfrak{G}_0$-measurable, and for any $m = 0, \ldots, M-1$, $\xi_m$ is both $\mathfrak{G}_{T_{m+1}}$-measurable and independent of $\mathfrak{G}_{T_m}$. Therefore, we only need to verify (16) for $X$ and $V$ defined on $(\bar{\Omega}, \mathfrak{G}, \bar{\mathbf{P}}, (\mathfrak{G}_t)_{t \geq 0})$.

Lemma 2 recalls well-known results. They may be deduced using Itô's formula, the existence of a smooth version of the stochastic flow $x \mapsto X_t^{s,x}$, and induction (some details may be found, e.g., in [18, 30]).

LEMMA 2. *Fix $\beta \in \mathbb{N} \cup \{0\}$. Suppose that $b$ and $\sigma$ are continuous functions such that $\partial_x^{\vec{p}} b$ and $\partial_x^{\vec{p}} \sigma$ are bounded continuous functions on $[0,T] \times \mathbb{R}^d$ for all $\vec{p} \in \mathcal{P}_l$, with $l = 1, \ldots, \beta + 2$. Let $(g_\theta)_{\theta \in \Theta}$ belong to $\mathcal{C}_p^{\beta+2}(\mathbb{R}^d, \mathbb{R})$. Set*

$$u_\theta(s, x) = \mathbf{E}(g_\theta(X_T) / X_s = x) = \mathbf{E}(g_\theta(X_T^{s,x})),$$

*whenever $s \in [0, T]$. Then $(u_\theta)_{\theta \in \Theta} \in \mathcal{C}_p^{\beta+2}(\mathbb{R}^d, \mathbb{R})$ and for all $\theta \in \Theta$,*

(17)
$$\frac{\partial}{\partial s} u_\theta(s, x) = -\mathcal{L}(u_\theta)(s, x) \quad \text{if } s \in [0, T] \text{ and } x \in \mathbb{R}^d,$$
$$u_\theta(T, x) = g_\theta(x) \quad \text{if } x \in \mathbb{R}^d,$$

*where*

$$\mathcal{L} = \sum_{k=1}^d b^k \, \partial_x^k + \tfrac{1}{2} \sum_{k,l=1}^d (\sigma \sigma^\top)^{k,l} \, \partial_x^{k,l}.$$

*Furthermore, for each $\vec{p} \in \mathcal{P}_l$, with $l = 0, \ldots, \beta$, we have $\frac{\partial}{\partial t} \partial_x^{\vec{p}} u_\theta$ is a continuous function and*

$$\frac{\partial}{\partial t} \partial_x^{\vec{p}} u_\theta = -\partial_x^{\vec{p}} \mathcal{L}(u_\theta).$$

In the sequel, the symbol $\mathbf{E}_t$ also denotes the conditional expectation with respect to $\mathfrak{G}_t$.

LEMMA 3. *Let $t \in [T_m, T_{m+1}]$. Then, for any $\mathfrak{G}_{T_m}$-measurable random variable $\xi$, we have*

(18) $\mathbf{E}_{T_m} g(t, Y_t^\xi) = g(T_m, \xi) + \mathbf{E}_{T_m} \int_{T_m}^t \left( \frac{\partial}{\partial t} g(s, Y_s^\xi) + \mathcal{L}_{T_m, \xi}(g)(s, Y_s^\xi) \right) ds,$

*with*

$$\mathcal{L}_{r,\xi}(g)(s,x) = \sum_{k=1}^{d} (Jb(r,\xi)x + b(r,\xi) - Jb(r,\xi)\xi)^k \, \partial_x^k g(s,x)$$

$$+ \tfrac{1}{2} \sum_{k,l=1}^{d} (\sigma\sigma^\top)^{k,l}(r,\xi) \, \partial_x^{k,l} g(s,x),$$

*provided that all of the derivates of* $g:[0,T] \times \mathbb{R}^d \to \mathbb{R}$ *appearing in* (18) *are continuous.*

PROOF. Observe that

$$Y_t^\xi = \xi + \int_{T_m}^t (Jb(T_m,\xi)Y_s^\xi + b(T_m,\xi) - Jb(T_m,\xi)\xi) \, ds$$

$$+ \int_{T_m}^t \sigma(T_m,\xi) \, dW_s.$$

Then we use the Itô formula to obtain (18). □

The proof of Lemma 4 is based on the Burkholder–Davis–Gundy inequalities and the discrete Gronwall–Bellman lemma. We omit it because it may be proved in much the same way as Lemma 4.3 of [24].

LEMMA 4. *Let assumptions of Lemma* 1 *hold. Then*

(19) $$\mathbf{E}\bigg(\sup_{m=0,\ldots,M} \|V_m^M\|^q\bigg) \leq K(T)(1 + \mathbf{E}(\|V_0^M\|^q))$$

*and*

(20) $\mathbf{E}_{T_m}(\|V_{m+1}^M - e^{Jb(T_m,V_m^M)T/M} V_m^M\|^q) \leq K(T)\Big(\dfrac{T}{M}\Big)^{q/2}(1 + \|V_m^M\|^q).$

*In addition, for any* $\mathfrak{G}_{T_m}$*-random variable* $\xi$ *taking values in* $\mathbb{R}^d$*, we have*

(21) $$\mathbf{E}_{T_m}\bigg(\sup_{t \in [T_m, T_{m+1}]} \|Y_t^{\xi,M} - e^{Jb(T_m,\xi)(t-T_m)}\xi\|^q\bigg) \leq K(T)\Big(\dfrac{T}{M}\Big)^{q/2}(1 + \|\xi\|^q).$$

The next result provides a uniform bound of weak order 1 for the weak error between $X_{T_m}$ and $V_m$, with $m = 0, \ldots, M$.

PROPOSITION 1. *Suppose that* $b$ *and* $\sigma$ *are continuous functions such that* $\partial_x^{\vec{p}} b$ *and* $\partial_x^{\vec{p}} \sigma$ *are bounded continuous functions on* $[0,T] \times \mathbb{R}^d$ *for any* $\vec{p} \in \mathcal{P}_l$*, with* $l = 1, \ldots, 4$*. Assume that the components of* $\partial b/\partial t$ *and* $\partial \sigma/\partial t$ belong to $\mathcal{C}_p^0([0,T] \times \mathbb{R}^d, \mathbb{R})$. Let $(g_\theta)_{\theta \in \Theta} \in \mathcal{C}_p^4(\mathbb{R}^d, \mathbb{R})$. If $V_0$ have the same distribution as $X_0$, then

$$|\mathbf{E}g_\theta(X_{T_m}) - \mathbf{E}g_\theta(V_m^M)| \leq K(T)(1 + \mathbf{E}\|X_0\|^q)\frac{T}{M}$$

for all $\theta \in \Theta$ and $m = 0, \ldots, M$.

PROOF. We first decompose the global error as a sum of terms involving the solution of a parabolic partial differential equation (17). This methodology goes back to Talay [27, 28] and Milshtein [21]. More precisely, let $u_\theta(s,x) = \mathbf{E}g_\theta(X_{T_m}^{s,x})$. Then

$$|\mathbf{E}g_\theta(X_{T_m}) - \mathbf{E}g_\theta(V_m)| \leq \sum_{j=1}^m |\mathbf{E}(H_1^j)| + |\mathbf{E}(H_2^j)|,$$

where $H_2^j = u_\theta(T_j, Y_{T_j}^{V_{j-1}}) - u_\theta(T_{j-1}, V_{j-1})$ and

$$H_1^j = u_\theta(T_j, V_j) - u_\theta(T_j, \tilde{V}_{j-1}) + u_\theta(T_j, \tilde{V}_{j-1}) - u_\theta(T_j, Y_{T_j}^{V_{j-1}}),$$

with $\tilde{V}_{j-1} = \exp(Jb(T_{j-1}, V_{j-1})T/M)V_{j-1}$.

Due to Lemma 2, we can use Taylor's formula to obtain

$$H_1^j = \sum_{l=1}^3 \frac{1}{l!} \sum_{\vec{p} \in \mathcal{P}_l} \partial_x^{\vec{p}} u_\theta(T_j, \tilde{V}_{j-1})(F_{\vec{p}}(V_j - \tilde{V}_{j-1}) - F_{\vec{p}}(Y_{T_j}^{V_{j-1}} - \tilde{V}_{j-1}))$$
$$+ R_j(V_j) + R_j(Y_{T_j}^{V_{j-1}}),$$

with

$$R_j(x) = \frac{1}{4!} \sum_{\vec{p} \in P_4} \partial_x^{\vec{p}} u_\theta(T_j, \tilde{V}_{j-1} + \Xi_{\vec{p},j}(x)(x - \tilde{V}_{j-1}))F_{\vec{p}}(x - \tilde{V}_{j-1}),$$

where $\Xi_{\vec{p},j}$ $d \times d$ is a diagonal matrix whose components belong to $[0,1]$. As in the proof of Theorem 4.1 of [24], using the Cauchy–Schwarz inequality and Lemmas 2 and 4, we now see that

$$(22) \qquad |\mathbf{E}H_1^j| \leq K(T)(\mathbf{E}\|V_{j-1}\|^q + 1)\left(\frac{T}{M}\right)^2.$$

Due to Lemma 2, applying Lemma 3 gives

$$\mathbf{E}_{T_{j-1}} H_2^j = \int_{T_{j-1}}^{T_j} \mathbf{E}_{T_{j-1}}(-\mathcal{L}(u_\theta)(s, Y_s^{V_{j-1}}) + \mathcal{L}_{T_{j-1}, V_{j-1}}(u_\theta)(s, Y_s^{V_{j-1}})) \, ds.$$

Then, combining again Lemmas 2 and 3 yields

$$\mathbf{E}_{T_{j-1}} H_2^j = \int_{T_{j-1}}^{T_j} \int_{T_{j-1}}^{t} \mathbf{E}_{T_{j-1}} \mathcal{L}_{T_{j-1}, V_{j-1}}^2(u_\theta)(s, Y_s^{V_{j-1}}) \, ds \, dt$$

$$+ \int_{T_{j-1}}^{T_j} \int_{T_{j-1}}^{t} \mathbf{E}_{T_{j-1}} (\mathcal{L}^2(u_\theta)(s, Y_s^{V_{j-1}})$$

(23)
$$- \mathcal{L}_1(u_\theta)(s, Y_s^{V_{j-1}})) \, ds \, dt$$

$$- 2 \int_{T_{j-1}}^{T_j} \int_{T_{j-1}}^{t} \mathbf{E}_{T_{j-1}} \mathcal{L}_{T_{j-1}, V_{j-1}}(\mathcal{L}(u_\theta))(s, Y_s^{V_{j-1}}) \, ds \, dt,$$

where

$$\mathcal{L}_1 = \sum_{k=1}^{d} \left( \frac{\partial}{\partial t} b^k \right) \partial_x^k + \frac{1}{2} \sum_{k,l=1}^{d} \left( \frac{\partial}{\partial t} (\sigma \sigma^\top)^{k,l} \right) \partial_x^{k,l}.$$

Hence, (21) and Lemma 2 lead to

(24) $$|\mathbf{E} H_2^j| \leq K(T)(\mathbf{E}\|V_{j-1}\|^q + 1) \left( \frac{T}{M} \right)^2.$$

From (19), (22) and (24) we deduce the assertion of this proposition. $\square$

We are now in position to show the main result of this section.

PROOF OF THEOREM 1. Let $u(s, x) = \mathbf{E} f(X_T^{s,x})$. As in the proof of Proposition 1, we have

$$\mathbf{E} f(V_M) = \mathbf{E} \sum_{m=1}^{M} (H_1^m + H_2^m) + \mathbf{E} u(0, V_0),$$

where $H_2^m = u(T_m, Y_{T_m}^{V_{m-1}}) - u(T_{m-1}, V_{m-1})$ and

$$H_1^m = u(T_m, V_m) - u(T_m, \tilde{V}_{m-1}) + u(T_m, \tilde{V}_{m-1}) - u(T_m, Y_{T_m}^{V_{m-1}}),$$

with $\tilde{V}_{m-1} = \exp(Jb(T_{m-1}, V_{m-1})T/M)V_{m-1}$.

It follows from Taylor's formula that

$$H_1^m = \sum_{l=1}^{5} \frac{1}{l!} \sum_{\vec{p} \in \mathcal{P}_l} \partial_x^{\vec{p}} u(T_m, \tilde{V}_{m-1})(F_{\vec{p}}(V_m - \tilde{V}_{m-1}) - F_{\vec{p}}(Y_{T_m}^{V_{m-1}} - \tilde{V}_{m-1}))$$

$$+ R_m^1(V_m) + R_m(Y_{T_m}^{V_{m-1}}),$$

with

$$R_m(x) = \frac{1}{6!} \sum_{\vec{p} \in P_6} \partial_x^{\vec{p}} u(T_m, \tilde{V}_{m-1} + \Xi_{\vec{p}, m}(x)(x - \tilde{V}_{m-1})) F_{\vec{p}}(x - \tilde{V}_{m-1}),$$

where $\Xi_{\vec{p},m}$ $d \times d$ is a diagonal matrix whose components belong to $[0,1]$. We may check that, for any $\vec{p} \in P_l$, with $l = 1, \ldots, 5$, there exist a constant $q$ and functions $c_{\vec{p}}$ in $C_p^4(\mathbb{R}^d, \mathbb{R})$ such that, for all $m = 1, \ldots, M$,

$$\left| \mathbf{E}_{T_{m-1}} F_{\vec{p}}(V_m - \tilde{V}_{m-1}) - \mathbf{E}_{T_k} F_{\vec{p}}(Y_{T_m}^{V_{m-1},M} - \tilde{V}_{m-1}) - c_{\vec{p}}(V_{m-1})\left(\frac{T}{M}\right)^2 \right|$$

$$\leq K(T)(1 + \|V_{m-1}\|^q)\left(\frac{T}{M}\right)^3.$$

Proceeding similarly to the proof of (22), we obtain

$$\left| \mathbf{E}\left( H_1^m - \left(\frac{T}{M}\right)^2 \sum_{l=1}^{5} \sum_{\vec{p} \in \mathcal{P}_l} \frac{\partial_x^{\vec{p}} u(T_m, \tilde{V}_{m-1})}{l!} c_{\vec{p}}(V_{m-1}) \right) \right|$$

$$\leq K(T)(\mathbf{E}\|V_{m-1}\|^q + 1)\left(\frac{T}{M}\right)^3.$$

Thus, the mean value theorem leads to

(25) $$\left| \mathbf{E}\left( H_1^m - \left(\frac{T}{M}\right)^2 \sum_{l=1}^{5} \sum_{\vec{p} \in \mathcal{P}_l} \frac{\partial_x^{\vec{p}} u(T_m, V_{m-1})}{l!} c_{\vec{p}}(V_{m-1}) \right) \right|$$

$$\leq K(T)(1 + \mathbf{E}\|V_{m-1}\|^q)\left(\frac{T}{M}\right)^3.$$

As in the estimation of $H_m^2$ in the proof of Proposition 1, Lemmas 2 and 3 imply that (23) holds with $u$ instead of $u_\theta$. Then, due to Lemma 2, applying Lemma 3 yields

$$\mathbf{E}_{T_{m-1}} H_2^m = \frac{T^2}{2M^2}\Lambda(T_{m-1}, V_{m-1}) + \mathbf{E}_{T_{m-1}} \int_{T_{m-1}}^{T_m} \int_{T_{m-1}}^{t} \int_{T_{m-1}}^{s} \Phi_m(r)\, dr\, ds\, dt,$$

where $\Lambda(s, \xi) = (\mathcal{L}_{s,\xi}^2(u) - \mathcal{L}_1(u) + \mathcal{L}^2(u) - 2\mathcal{L}_{s,\xi}(\mathcal{L}(u)))(s, \xi)$ and

$$\Phi_m(s) = (2\mathcal{L}_1(\mathcal{L}(u)) + \mathcal{L}(\mathcal{L}_1(u)) - \mathcal{L}^3(u) - \mathcal{L}_2(u))(s, Y_s^{V_{m-1}})$$
$$+ (3\mathcal{L}_{T_{m-1}, Z_{m-1}}(\mathcal{L}^2(u)) - 3\mathcal{L}_{T_{m-1}, V_{m-1}}(\mathcal{L}_1(u)))(s, Y_s^{V_{m-1}})$$
$$+ (\mathcal{L}_{T_{m-1}, V_{m-1}}^3(u) - 3\mathcal{L}_{T_{m-1}, V_{m-1}}^2(\mathcal{L}(u)))(s, Y_s^{V_{m-1}}).$$

Here

$$\mathcal{L}_2 = \sum_{k=1}^{d} \left(\frac{\partial^2}{\partial t^2} b^k\right) \partial_x^k + \frac{1}{2} \sum_{k,l=1}^{d} \left(\frac{\partial^2}{\partial t^2}(\sigma\sigma^\top)^{k,l}\right) \partial_x^{k,l}.$$

It follows from Lemma 2 that

(26) $$\left| \mathbf{E}\left( H_2^m - \left(\frac{T}{M}\right)^2 \Lambda(T_{m-1}, V_{m-1})/2 \right) \right| \leq K(T)(\mathbf{E}\|V_{m-1}\|^q + 1)\left(\frac{T}{M}\right)^3.$$

Let

$$\Psi(s,x) = \Lambda(s,x)/2 + \sum_{l=1}^{5} \sum_{\vec{p} \in \mathcal{P}_l} \partial_x^{\vec{p}} u(s,x) c_{\vec{p}}(x)/l!.$$

Using Lemma 2, we get that $(\Psi(s,\cdot))_{s \in [0,T]} \in C_p^4(\mathbb{R}^d, \mathbb{R})$. Then, it follows from Proposition 1 that

$$(27) \quad |\mathbf{E}(\Psi(T_{m-1}, V_{m-1}) - \Psi(T_{m-1}, X_{T_{m-1}}))| \leq K(T)(\mathbf{E}\|X_0\|^q + 1)\frac{T}{M}.$$

Hence, combining Itô's formula, (25), (26) and (27) give

$$\left| \mathbf{E}\left( H_1^m + H_2^m - \frac{T}{M} \int_{T_{m-1}}^{T_m} \Psi(s, X_s) \, ds \right) \right| \leq K(T)(\mathbf{E}\|X_0\|^q + 1)\left(\frac{T}{M}\right)^3,$$

which completes the proof. □

REMARK 3. According to Theorem 1, we have

$$|\mathbf{E}f(X_T) - 2\mathbf{E}f(V_{2M}^{2M}) + \mathbf{E}f(V_M^M)| \leq K(T)(1 + \mathbf{E}\|X_0\|^q)\left(\frac{T}{M}\right)^2,$$

provided the hypotheses of Theorem 1 hold. This yields a second weak order scheme based on the extrapolation Scheme 1 (see, e.g., [32]).

**4. Euler-exponential scheme for stochastic Schrödinger equations.** We now turn to our main problem. To be more precise, this section provides an heuristic deduction of a version of the Euler-exponential scheme adapted to the characteristics of (1).

The following lemma discusses the existence and uniqueness of the solutions of (1).

LEMMA 5. *Let $s \geq 0$. Suppose that $\xi$ is a $\mathfrak{F}_s$-random variable with $\mathbf{E}\|\xi\|^2 < \infty$. Then there exists a unique global continuous solution of the SDE*

$$(28) \quad Z_t^{s,\xi} = \xi + \int_s^t (GZ_r^{s,\xi} + D(Z_r^{s,\xi})) \, dr + \sum_{k=1}^n \int_s^t E_k(Z_r^{s,\xi}) \, dW_r^k$$

*for all $t \geq s$. Moreover, $\|Z^{s,\xi}\| = 1$ a.s., provided that $\|\xi\| = 1$ a.s. Recall that $D$ is given by (2) and $E_1, \ldots, E_n$ are defined by (3).*

PROOF. Since the drift coefficient of (28) and $E_k$, with $k = 1, \ldots, n$, are locally Lipschitz, applying the truncation method, we obtain that (28) has a unique local solution (see, e.g., [19, 26]). That is, there exists a stopping time $\zeta_\xi$ such that (28) has a unique solution up to $\zeta_\xi$. This solution has continuous paths a.s. and $\limsup_{t \to \zeta_\xi-} \|Z_t^{s,\xi}\| = \infty$ a.s. on $\{\zeta_\xi < \infty\}$. By Itô's formula, we have

$$(29) \quad \|Z_{\tau_N \wedge t}^{s,\xi}\|^2 = \|\xi\|^2 + 2 \sum_{k=1}^{n} \int_{s}^{\tau_N \wedge t} \operatorname{Re} \langle Z_r^{s,\xi}, L_k Z_r^{s,\xi} \rangle (1 - \|Z_r^{s,\xi}\|^2) \, dW_r^k,$$

where $\tau_N$ is the first exit time of $Z_t^{s,\xi}$ of $\{x : \|x\| \le N\}$. This yields

$$\mathbf{E} \|Z_{\tau_N \wedge t}^{s,\xi}\|^2 = \mathbf{E} \|\xi\|^2.$$

It follows that $\zeta_\xi = +\infty$ a.s. (null set depending on $\xi$). Hence, there exists a unique global continuous solution of (28).

Let $S_t = \sum_{k=1}^{n} \int_s^t \operatorname{Re} \langle Z_r^{s,\xi}, L_k Z_r^{s,\xi} \rangle \, dW_r^k$. Then $(S_t)_{t \ge s}$ is a continuous semimartingale and

$$(30) \quad \|Z_t^{s,\xi}\|^2 = \|\xi\|^2 + \int_s^t (1 - \|Z_r^{s,\xi}\|^2) \, dS_r.$$

Thus, the last assertion of the lemma follows from the uniqueness of the solution of (30). □

We split the drift coefficient of (1) into $GZ_s$ and $D(Z_s)$. Then, analysis similar to that in the proof of Lemma 1 shows that, for all $t \in [T_m, T_{m+1}]$,

$$
\begin{aligned}
Z_t^{T_m, Z_{T_m}} &= \exp(G(t - T_m)) Z_{T_m} + \int_{T_m}^t \exp(G(t-s)) D(Z_s^{T_m, Z_{T_m}}) \, ds \\
&\quad + \sum_{k=1}^n \int_{T_m}^t \exp(G(t-s)) E_k(Z_s^{T_m}, Z_{T_m}) \, dW_s^k.
\end{aligned}
$$
(31)

Let $\hat{Z}_m$ be a linear weak approximation of $Z_{T_m}$ satisfying $\|\hat{Z}_m\| = 1$. We now approximate $Z_s^{T_m, Z_{T_m}}$ by $\hat{Z}_m$ in the right-hand side of (31) to obtain

$$Z_t \approx e^{G(t-T_m)} \hat{Z}_m + \int_{T_m}^t e^{G(t-s)} D(\hat{Z}_m) \, ds + \sum_{k=1}^n \int_{T_m}^t e^{G(t-s)} E_k(\hat{Z}_m) \, dW_s^k$$
(32)
for all $t \in [T_m, T_{m+1}]$.

Since our goal is to compute $\mathbf{E}\langle Z_t, A Z_t \rangle$, we should approximate the measure induced by the right-hand side of (32). To this end, we use the procedure employed in Section 3.1 to yield

$$Z_{T_{m+1}} \approx \bar{Z}_{m+1} = \Phi_{m+1}^{\hat{Z}_m, M},$$

where $\Phi$ is given by (7). Finally, to include the information that $\|Z_t\| = 1$, $\bar{Z}_{m+1}$ is projected onto the manifold $\{z \in \mathbb{C}^d : \|z\| = 1\}$. We thus get Scheme 2 defined in the Introduction, that is, the method

$$\hat{Z}_{m+1}^M = p \left( \exp\left( G \frac{T}{M} \right) \left( \hat{Z}_m^M + D(\hat{Z}_m^M) \frac{T}{M} + \sqrt{\frac{T}{M}} \sum_{k=1}^n E_k(\hat{Z}_m^M) \xi_m^k \right) \right),$$

where

$$p(z) = \begin{cases} 0, & \text{if } z = 0, \\ z/\|z\|, & \text{if } z \neq 0. \end{cases}$$

REMARK 4. As in Remark 1, Scheme 2 leads to the following version of the implicit Euler scheme.

SCHEME 3. Let $\hat{I}_0^M$ be a random variable with $\|\hat{I}_0^M\| = 1$. Then we set

$$\hat{I}_{m+1}^M = p\left(\left(I - G\frac{T}{M}\right)^{-1}\left(\hat{I}_m^M + D(\hat{I}_m^M)\frac{T}{M} + \sqrt{\frac{T}{M}}\sum_{k=1}^n E_k(\hat{I}_m^M)\xi_m^k\right)\right).$$

**5. Rate of convergence.** In this section we focus our interest on the proof of the following theorem which establishes the linear convergence of $\mathbf{E}\langle \hat{Z}_M^M, A\hat{Z}_M^M \rangle$.

THEOREM 2. *Suppose that, for all $B \in \mathbb{C}^{d,d}$,*

$$(33) \qquad |\mathbf{E}\langle Z_0, BZ_0 \rangle - \mathbf{E}\langle \hat{Z}_0^M, B\hat{Z}_0^M \rangle| \leq \|B\|K(T)\frac{T}{M}.$$

*If the law of $\xi_0^1$ has compact support, then*

$$(34) \qquad |\mathbf{E}\langle Z_T, AZ_T \rangle - \mathbf{E}\langle \hat{Z}_M^M, A\hat{Z}_M^M \rangle| \leq K(T)\frac{T}{M}.$$

The proof of Theorem 2 starts with bounds for the concentration of $\Phi$. As in [9, 29], the assumption that $\xi_0^1$ has compact support allows us to obtain this kind of estimate.

LEMMA 6. *Let $\xi_0^1$ have compact support. Then there exists an increasing positive function $K_2$ such that $\|\Phi_{m+1}^{z,M}\| \leq K_2(T)$, whenever $m = 0, \ldots, M-1$ and $\|z\| = 1$. Furthermore, there exist $\delta \in {]0,1[}$ and a strictly positive constant $K_1$ independent of both $T$ and $T/M$ such that $\|\Phi_{m+1}^{z,M}\| \geq K_1$ for all $T/M < \delta$, $m = 0, \ldots, M-1$, and $z \in \mathbb{C}^d$ with $\|z\| = 1$.*

PROOF. Without loss of generality, we can assume that the support of $\xi_0^1$ belongs to the interval $[-a, a]$. Hence, for any $z \in \mathbb{C}^d$, with $\|z\| = 1$, we get

$$(35) \qquad \left\|\frac{T}{M}D(z) + \sqrt{\frac{T}{M}}\sum_{k=1}^n E_k(z)\xi_m^k\right\| \leq \frac{3}{2}\frac{T}{M}\sum_{k=1}^n \|L_k\|^2 + 2a\sqrt{\frac{T}{M}}\sum_{k=1}^n \|L_k\|.$$

For any $z \in \mathbb{C}^d$, $\operatorname{Re}\langle Gz, z \rangle \leq 0$, and so $(\exp(Gt))_{t \geq 0}$ is a contraction semigroup on $\mathbb{C}^d$. This yields

$$\|\Phi_{m+1}^{z,M}\| \leq \left\| z + D(z)\frac{T}{M} + \sqrt{\frac{T}{M}} \sum_{k=1}^n E_k(z)\xi_m^k \right\| \leq K_2(T),$$

for all $z \in \mathbb{C}^d$ with $\|z\| = 1$.

Let $\phi(z) = \|\exp(G)z\|$. From (35), we see that there exists $\delta \in\,]0,1[$ such that

$$\left\| z + D(z)\frac{T}{M} + \sqrt{\frac{T}{M}} \sum_{k=1}^n E_k(z)\xi_m^k \right\| \geq \frac{1}{2},$$

for $T/M < \delta$ and $z \in \mathbb{C}^d$ with $\|z\| = 1$. Thus,

$$\phi\left( z + D(z)\frac{T}{M} + \sqrt{\frac{T}{M}} \sum_{k=1}^n E_k(z)\xi_m^k \right) \geq K_1,$$

provided that $\|z\| = 1$ and $T/M < \delta$. This gives the last assertion of the lemma since $\phi(x) \leq \|\exp(GT/M)x\|$ whenever $T/M < 1$. $\square$

Similarly to the proof of Theorem 1, in the sequel we consider the complete probability space $(\bar{\Omega}, \mathfrak{G}, \bar{\mathbf{P}})$ induced by the random variables $\hat{Z}_0^M$, $\xi_0, \ldots, \xi_{M-1}$. In addition, $(\bar{\Omega}, \mathfrak{G}, \bar{\mathbf{P}}, (\mathfrak{G}_t)_{t \geq 0})$ will be the filtered probability space satisfying the usual hypotheses induced by $(\bar{\Omega}, \mathfrak{G}, \bar{\mathbf{P}})$ and the filtration $(\sigma(\hat{Z}_0^M, \xi_k : k \leq [tM/T] - 1))_{t \geq 0}$. By abuse of notation, we use the same symbol $\mathbf{E}_t$ for the conditional expectation with respect to both $\mathfrak{F}_t$ and $\mathfrak{G}_t$.

The role of the local approximation $Y$ in the proof of Theorem 1 is played here by $\Psi$ given by

$$\Psi_t^{s,z} = z + \int_s^t (G\Psi_{r_-}^{s,z} + D(z))\,dr + \sum_{k=1}^n \int_s^t E_k(z)\,dS_r^k,$$

where $t \geq s$ and for any $k = 1, \ldots, n$,

$$S_r^k = \sqrt{\frac{T}{M}} \sum_{m=0}^{M-1} \xi_m^k I_{[T_{m+1}, +\infty[}(r).$$

The next two lemmas provide information about the behavior of $\Psi$.

LEMMA 7. *Let $\xi_0^1$ have compact support. Then there exists an increasing positive function $K_4$ such that $\|\Psi_t^{T_m^M, z}\| \leq K_4(T)$, whenever $\|z\| = 1$, $m = 0, \ldots, M-1$, and $t \in [T_m^M, T_{m+1}^M]$. Moreover, there exist $\Delta \in\,]0,1[$ and a strictly positive constant $K_3$ independent of both $T$ and $T/M$ such that $\|\Psi_t^{T_m^M, z}\| \geq K_3$ for all $T/M < \Delta$, $m = 0, \ldots, M-1$, $t \in [T_m^M, T_{m+1}^M]$, and $z \in \mathbb{C}^d$ with $\|z\| = 1$.*

PROOF. For any $t \in [T_m, T_{m+1}]$, we have

$$\Psi_t^{T_m,z} = e^{G(t-T_m)}\left(z + \int_{T_m}^t e^{G(T_m-s)} D(z)\, ds \right.$$

$$\left. + \sqrt{\frac{T}{M}} \sum_{k=1}^n e^{G(T_m-t)} E_k(z) \xi_m^k I_{\{T_{m+1}\}}(t)\right).$$

Hence, analysis similar to that in the proof of Lemma 6 shows the assertion of the lemma. $\square$

LEMMA 8. *Suppose that $K_3, K_4$ and $\Delta$ are as in Lemma 7. Assume that $f \in C^{2,4,4}([0,T] \times \mathcal{S} \times \mathcal{S}, \mathbb{C})$, with $\mathcal{S} = \{z \in \mathbb{C}^d : K_3(T) \leq \|z\| \leq K_4(T)\}$. Let $\xi_0^1$ have compact support. Then for all $t \in [T_m, T_{m+1}[$,*

$$\mathbf{E}_{T_m} f(t, \Psi_t^{T_m,z}, \overline{\Psi_t^{T_m,z}})$$
(36)
$$= f(T_m, z, \bar{z})$$
$$+ \int_{T_m}^t \mathbf{E}_{T_m}\left(\frac{\partial f}{\partial s}(s, \Psi_s^{T_m,z}, \overline{\Psi_s^{T_m,z}}) + \mathcal{L}_z^1(f)(s, \Psi_s^{T_m,z}, \overline{\Psi_s^{T_m,z}})\right) ds,$$

*provided that $T/M < \Delta$ and $z \in \mathbb{C}^d$ with $\|z\| = 1$. Here*

$$\mathcal{L}_z^1(f)(s,x,y) = \sum_{k=1}^d \left(\frac{\partial f}{\partial x^k}(s,x,y)(Gx + D(z))^k + \frac{\partial f}{\partial y^k}(s,x,y)(\bar{G}y + \overline{D(z)})^k\right).$$

*Furthermore,*

$$\mathbf{E}_{T_m} f(T_{m+1}, \Psi_{T_{m+1}}^{T_m,z}, \overline{\Psi_{T_{m+1}}^{T_m,z}})$$
$$= f(T_m, z, \bar{z})$$
(37)
$$+ \int_{T_m}^{T_{m+1}} \mathbf{E}_{T_m}\left(\frac{\partial f}{\partial s}(s, \Psi_s^{T_m,z}, \overline{\Psi_s^{T_m,z}}) + \mathcal{L}_z(f)(s, \Psi_s^{T_m,z}, \overline{\Psi_s^{T_m,z}})\right) ds$$
$$+ O_f\left(z, \frac{T}{M}\right),$$

*where $\|O_f(z, T/M)\| \leq K(T)(T/M)^2$ and*

$$\mathcal{L}_z(f)(s,x,y) = \mathcal{L}_z^1(f)(s,x,y) + \frac{1}{2}\sum_{j=1}^n \sum_{k,l=1}^d \frac{\partial^2 f}{\partial x^k x^l}(s,x,y) E_j(z)^k E_j(z)^l$$
$$+ \sum_{j=1}^n \sum_{k,l=1}^d \left(\frac{\partial^2 f}{\partial x^k y^l}(s,x,y) E_j(z)^k \overline{E_j(z)}^l\right.$$
$$\left. + \frac{1}{2}\frac{\partial^2 f}{\partial y^k y^l}(s,x,y) \overline{E_j(z)}^k \overline{E_j(z)}^l\right).$$

PROOF. Applying the Itô formula for a general semimartingale, we obtain, for any $t \in [T_m, T_{m+1}]$,

$$\mathbf{E}_{T_m} f(t, \Psi_t^{T_m,z}, \overline{\Psi_t^{T_m,z}}) - f(T_m, z, \bar{z})$$
$$= \int_{T_m}^t \mathbf{E}_{T_m} \left( \frac{\partial f}{\partial s}(s, \Psi_s^{T_m,z}, \overline{\Psi_s^{T_m,z}}) + \mathcal{L}_z^1(f)(s, \Psi_s^{T_m,z}, \overline{\Psi_s^{T_m,z}}) \right) ds$$
$$+ I_{\{T_{m+1}\}}(t) \mathbf{E}_{T_m} \Bigg( f(t, \Psi_t^{T_m,z}, \overline{\Psi_t^{T_m,z}}) - f(t, \Psi_{t-}^{T_m,z}, \overline{\Psi_{t-}^{T_m,z}})$$
$$- \sum_{k=1}^d \frac{\partial f}{\partial x^k}(t, \Psi_{t-}^{T_m,z}, \overline{\Psi_{t-}^{T_m,z}})(\Psi_t^{T_m,z} - \Psi_{t-}^{T_m,z})^k$$
$$- \sum_{k=1}^d \frac{\partial f}{\partial y^k}(t, \Psi_{t-}^{T_m,z}, \overline{\Psi_{t-}^{T_m,z}})(\overline{\Psi_t^{T_m,z}} - \overline{\Psi_{t-}^{T_m,z}})^k \Bigg).$$

This gives (36). Furthermore, expanding

$$f(T_{m+1}, \Psi_{T_{m+1}}^{T_m,z}, \overline{\Psi_{T_{m+1}}^{T_m,z}}) - f(T_{m+1}, \Psi_{T_{m+1}-}^{T_m,z}, \overline{\Psi_{T_{m+1}-}^{T_m,z}})$$

in powers of $(\Psi_{T_{m+1}}^{T_m,z} - \Psi_{T_{m+1}-}^{T_m,z})^j$ and $(\overline{\Psi_{T_{m+1}}^{T_m,z}} - \overline{\Psi_{T_{m+1}-}^{T_m,z}})^j$, with $j = 1, \ldots, d$, we obtain (37). To this end, we combine (36), Taylor's formula and the mean value theorem. □

Note that the coefficients of (1) are not globally Lipschitz. To overcome this difficulty, instead of using the solution of the usual (in this context) partial differential equation associated to (28), we employ the function $v \colon [0, T] \times \mathcal{D} \to \mathbb{C}$ described by $\mathcal{D} = \{(x, y) : x, y \in \mathbb{C}^d, \langle \bar{y}, x \rangle \neq 0\}$ and

$$v(s, x, y) = \langle \bar{y}, \tau_{T-s} x \rangle / \langle \bar{y}, x \rangle,$$

where $\tau_t$ is the solution of the backward quantum master equation (8) with $\mathfrak{h} = \mathbb{C}^d$.

PROOF OF THEOREM 2. Let $\alpha \colon [0, T] \times \mathbb{C}^d \times \mathbb{C}^d \mapsto \mathbb{C}$ be given by

$$\alpha(t, x, y) = \langle \bar{y}, \tau_{T-t} x \rangle.$$

According to (8) and Itô's formula, we have $\mathbf{E}\alpha(t, Z_t^{s,z}, \overline{Z_t^{s,z}}) = \alpha(s, z, \bar{z})$, for any $t \in [0, T]$. Hence, $\mathbf{E}\langle Z_T^{s,z}, A Z_T^{s,z} \rangle = \langle z, \tau_{T-s} z \rangle$. We thus get

$$v(s, z, \bar{z}) = \mathbf{E}\langle Z_T^{s,z}, A Z_T^{s,z} \rangle,$$

provided that $\|z\| = 1$. Therefore,

(38) $\quad \mathbf{E}\langle \hat{Z}_M^M, A \hat{Z}_M^M \rangle - \mathbf{E}\langle Z_T, A Z_T \rangle = \mathbf{E} v(T, \hat{Z}_M^M, \overline{\hat{Z}_M^M}) - \mathbf{E} v(0, Z_0, \overline{Z_0}).$

Suppose that $T/M < min\{\delta, \Delta\}$, with $\delta$ and $\Delta$ as in Lemma 6 and Lemma 7, respectively. Since $\|\hat{Z}_0\| = 1$, from (8), (33) and (38), we conclude that

$$\left\| \mathbf{E}\langle \hat{Z}_M^M, A\hat{Z}_M^M \rangle - \mathbf{E}\langle Z_T, AZ_T \rangle - \sum_{m=1}^{M} \mathbf{E}(H_1^m + H_2^m) \right\| \leq K(T)\frac{T}{M},$$

where

$$H_1^m = v\Big(T_m, \Phi_m^{\hat{Z}_{m-1},M}, \overline{\Phi_m^{\hat{Z}_{m-1},M}}\Big) - v(T_m, e^{GT/M}\hat{Z}_{m-1}, \overline{e^{GT/M}\hat{Z}_{m-1}})$$
$$+ v(T_m, e^{GT/M}\hat{Z}_{m-1}, \overline{e^{GT/M}\hat{Z}_{m-1}}) - v\Big(T_m, \Psi_{T_m}^{T_{m-1},\hat{Z}_{m-1}}, \overline{\Psi_{T_m}^{T_{m-1},\hat{Z}_{m-1}}}\Big)$$

and

$$H_2^m = v\Big(T_m, \Psi_{T_m}^{T_{m-1},\hat{Z}_{m-1}}, \overline{\Psi_{T_m}^{T_{m-1},\hat{Z}_{m-1}}}\Big) - v(T_{m-1}, \hat{Z}_{m-1}, \overline{\hat{Z}_{m-1}}).$$

We proceed to estimate $H_1^m$. From the construction of Scheme 2 and a simple computation, we see that, for any $z \in \mathbb{C}^d$ satisfying $\|z\| = 1$ and $\vec{p} \in \mathcal{P}_l$ with $l = 1, 2, 3$, we have

$$|\mathbf{E}_{T_{m-1}} F_{\vec{p}}(\Phi_m^{z,M} - e^{GT/M}z) - \mathbf{E}_{T_{m-1}} F_{\vec{p}}(\Psi_{T_m}^{T_{m-1},z} - e^{GT/M}z)| \leq K(T)\left(\frac{T}{M}\right)^2.$$

Hence, combining Lemma 6 with the deterministic Taylor formula gives

$$(39) \qquad |\mathbf{E} H_1^m| \leq K(T)\left(\frac{T}{M}\right)^2.$$

This follows by the same method as in the estimation of $H_1^m$ in the proof of Proposition 1.

It remains to estimate $H_2^m$. According to Lemma 8, we have

$$(40) \qquad \begin{aligned} \mathbf{E}_{T_{m-1}} H_2^m &= \mathbf{E}_{T_{m-1}} \int_{T_{m-1}}^{T_m} \frac{\partial v}{\partial s}\Big(s, \Psi_s^{T_m,\hat{Z}_{m-1}}, \overline{\Psi_s^{T_m,\hat{Z}_{m-1}}}\Big) ds \\ &\quad + \mathbf{E}_{T_{m-1}} \int_{T_{m-1}}^{T_m} \mathcal{L}_{\hat{Z}_{m-1}}(v)\Big(s, \Psi_s^{T_m,\hat{Z}_{m-1}}, \overline{\Psi_s^{T_m,\hat{Z}_{m-1}}}\Big) ds \\ &\quad + O_v\Big(\hat{Z}_{m-1}, \frac{T}{M}\Big). \end{aligned}$$

We may now apply Lemma 8 to the terms of the right-hand side of (40) to obtain

$$(41) \qquad |\mathbf{E} H_2^k| \leq K(T)\left(\frac{T}{M}\right)^2.$$

To be more precise, fortunately a very long computation shows that

$$\mathcal{L}_z(v)(s, z, \bar{z}) = \langle z, \tau_{T-s} Gz \rangle + \langle Gz, \tau_{T-s} z \rangle + \sum_{k=1}^{n} \langle L_k z, \tau_{T-s} L_k z \rangle,$$

whenever $\|z\| = 1$. Since (8) leads to

(42)
$$\frac{\partial v}{\partial t}(s,x,y) = -\frac{1}{\langle \bar{y}, x \rangle}\left(\langle \bar{y}, G^*\tau_{T-s}x\rangle + \langle \bar{y}, \tau_{T-s}Gx\rangle + \sum_{j=1}^{n}\langle \bar{y}, L_j^*\tau_{T-s}L_jx\rangle\right),$$

we deduce that, for $\|z\| = 1$,

$$\frac{\partial v}{\partial t}(s,z,\bar{z}) + \mathcal{L}_z(v)(s,z,\bar{z}) = 0.$$

Therefore, Lemma 8 yields that $\mathbf{E}_{T_{m-1}}H_2^m$ is equal to

$$O_v\left(\hat{Z}_{m-1}, \frac{T}{M}\right) + \mathbf{E}_{T_{m-1}}\int_{T_{m-1}}^{T_m}\int_{T_{m-1}}^{s}\frac{\partial^2 v}{\partial r}\left(r, \Psi_r^{T_m,\hat{Z}_{m-1}}, \overline{\Psi_r^{T_m,\hat{Z}_{m-1}}}\right)dr\,ds$$

$$+ \mathbf{E}_{T_{m-1}}\int_{T_{m-1}}^{T_m}\int_{T_{m-1}}^{s}\mathcal{L}_{\hat{Z}_{m-1}}^1\left(\frac{\partial v}{\partial r}\right)\left(r, \Psi_r^{T_m,\hat{Z}_{m-1}}, \overline{\Psi_r^{T_m,\hat{Z}_{m-1}}}\right)dr\,ds$$

$$+ \mathbf{E}_{T_{m-1}}\int_{T_{m-1}}^{T_m}\int_{T_{m-1}}^{s}\frac{\partial}{\partial r}(\mathcal{L}_{\hat{Z}_{m-1}}(v))\left(r, \Psi_r^{T_m,\hat{Z}_{m-1}}, \overline{\Psi_r^{T_m,\hat{Z}_{m-1}}}\right)dr\,ds$$

$$+ \mathbf{E}_{T_{m-1}}\int_{T_{m-1}}^{T_m}\int_{T_{m-1}}^{s}\mathcal{L}_{\hat{Z}_{m-1}}^1(\mathcal{L}_{\hat{Z}_{m-1}}(v))\left(r, \Psi_r^{T_m,\hat{Z}_{m-1}}, \overline{\Psi_r^{T_m,\hat{Z}_{m-1}}}\right)dr\,ds.$$

Hence, (41) follows from (42) and Lemma 7.

We conclude from (39) and (41) that (34) holds for $T/M < min\{\delta, \Delta\}$. Hence, our claim follows from $\|Z_T\| = 1$ a.s. and $\|\hat{Z}_M^M\| = 1$. □

REMARK 5. We expect that an expansion similar to (16) holds for the error $\mathbf{E}\langle Z_T, AZ_T\rangle - \mathbf{E}\langle \hat{Z}_M^M, A\hat{Z}_M^M\rangle$. Nevertheless, the proof of this result is still in progress.

REMARK 6. We now turn to (9) with $\dim \mathfrak{h} = +\infty$. It is relevant to characterize the global error $|\mathbf{E}_\mathbf{Q}\langle Y_T, AY_T\rangle - \mathbf{E}\langle \hat{Z}_M^M, A\hat{Z}_M^M\rangle|$ in function of $M$ and the dimension of $\mathfrak{h}_m$. A step toward this goal was given in [23], where the rate of convergence of $\mathbf{E}\langle E_{T,m}^M, AE_{T,m}^M\rangle$ to $\mathbf{E}_\mathbf{Q}\langle Y_T, AY_T\rangle$ is studied. Here $E$ denotes the numerical solution of (9) by the Euler scheme. An objective of this paper is to advance toward the solution of this problem.

**6. Numerical experiment.** This section illustrates the performance of the scheme $\hat{Z}$. To this end, we consider the following representative example of forced and damped quantum harmonic oscillator in the interaction representation.

EXAMPLE 1. Returning to Section 1.2, we choose $\mathfrak{h} = l^2(\mathbb{Z}_+)$. Let $(\varphi_k)_{k \in \mathbb{Z}_+}$ be the canonical orthonormal basis on the space $l^2(\mathbb{Z}_+)$. Then, the domain of the operators $a^\dagger$ and $a$ is $\{x \in l^2(\mathbb{Z}_+) : \sum_{k \geq 0} k|x_k|^2 < +\infty\}$ and for all $m \in \mathbb{Z}_+$, $a^\dagger \varphi_m = \sqrt{m+1} \varphi_{m+1}$ and

$$a\varphi_m = \begin{cases} 0, & \text{if } m = 0, \\ \sqrt{m}\varphi_{m-1}, & \text{if } m > 0. \end{cases}$$

The Number operator is defined by $N = a^\dagger a$.

We now simulate the Hamiltonian as $H = i(a^\dagger - a) + N$. Furthermore, we set $L_1 = 0.2a$, $L_2 = 0.01a^2$, $L_3 = 0.1N$ and $L_4 = 0.1a^\dagger$.

In Example 1, $\mathfrak{h}$ describes the state space of a single mode of a quantized electromagnetic field. The operator $a^\dagger$, respectively $a$, is the creation operator, respectively annihilation operator. Then, for instance, the term $i(a^\dagger - a)$ describes a linear pumping and $L_1$ simulates the damping due to photon emission.

To test the scheme $\hat{Z}$, we set $T = 100$ and $Y_0 = \varphi_6$. Moreover, we choose $\mathfrak{h}_d$ as the linear manifold spanned by $\{\varphi_j : 0 \leq j \leq d\}$ with $d = 50$. The objective is to describe numerically $\mathbf{E}_\mathbf{Q}\langle Y_{t,50}, NY_{t,50}\rangle$ for $t \in [0,T]$. As we comment in Section 1.2, this task leads to solve (1) for $d = 50$. The parameters selected allow us to obtain the "true" value of $\mathbf{E}_\mathbf{Q}\langle Y_{t,50}, NY_{t,50}\rangle$ by means of the solution of the backward quantum master equation (8) associated to our SDE. To this end, we use its explicit solution. It is worth pointing out that the numerical solution of finite-dimensional backward quantum master equations presents serious drawbacks when the dimension of the state space is high (see, e.g., [25]). In fact, some of these problems can be observed in our example in case $d = 100$.

In the numerical experiment, we compare Scheme 2, Scheme 3 and the following version of the explicit Euler scheme:

$$\hat{E}_{k+1} = \begin{cases} 0, & \text{if } \bar{E}_{k+1} = 0, \\ \bar{E}_{k+1}/\|\bar{E}_{k+1}\|, & \text{if } \bar{E}_{k+1} \neq 0, \end{cases}$$

where $\bar{E}_{k+1} = \widehat{E}_k + (G\widehat{E}_k + D(\widehat{E}_k))T/M + \sqrt{T/M} \sum_{j=1}^n E_j(\widehat{E}_k)\xi_k^j$, with $\xi_0^1, \ldots, \xi_0^n, \ldots, \xi_{M-1}^1, \ldots, \xi_{M-1}^n$ as in Scheme 2. In all codes, $\xi_0^1$ assumes values $\pm 1$, each with probability $1/2$.

Figure 1 shows the "true" solution and the approximations obtained by the numerical schemes. Moreover, Table 1 looks at the dependence of the errors $\varepsilon_0$ to the time step size $T/M$, where

$$\varepsilon_J(\chi, M) = \max_{j=J,\ldots,100} \left| \mathbf{E}\langle Z_j, NZ_j\rangle - 2\cdot 10^{-4} \sum_{k=1}^{2\cdot 10^{-4}} \langle \chi_{jM/100}^M(\omega_k), N\chi_{jM/100}^M(\omega_k)\rangle \right|,$$

NUMERICAL SOLUTION OF SSES 25

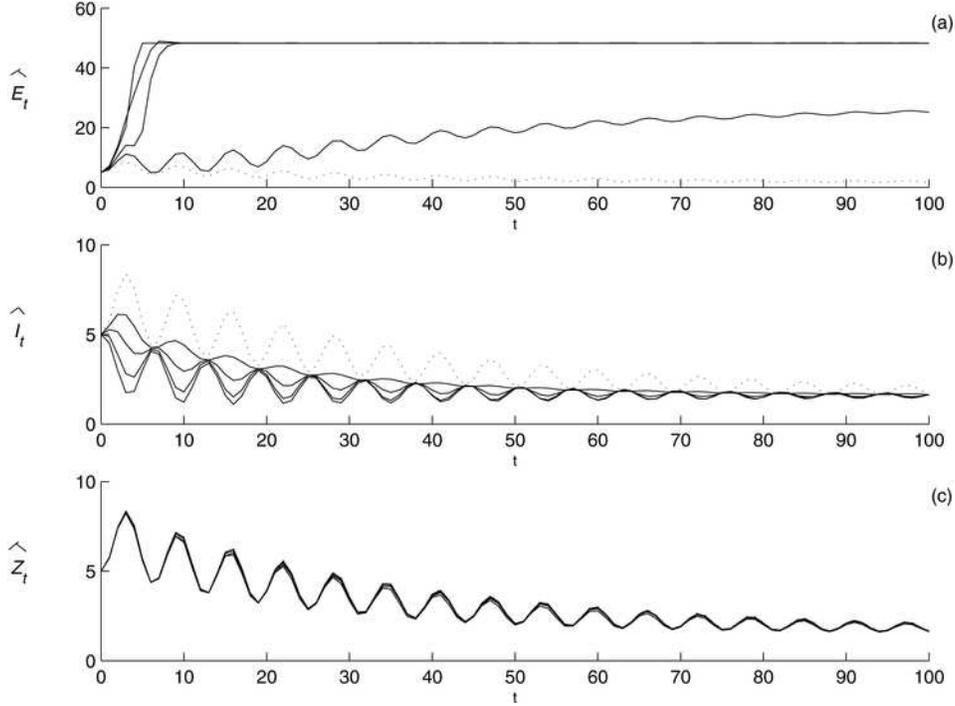

FIG. 1. *Dotted line:* "*true*" *solution, solid line:* (a) *explicit Euler scheme,* (b) *implicit Euler scheme and* (c) *Euler-exponential scheme.*

whenever $\chi$ denotes the numerical method and $J \in \{0, \ldots, 100\}$. Indeed, Table 1 presents estimated values of $\varepsilon_0$ and $\Delta$. Here $\Delta$ is the maximum of the length of the 90 percent confidence intervals taken over the instant of times $\{0, \ldots, 100\}$. We use the batch method to estimate these intervals (see, e.g., [17]).

TABLE 1
*Errors versus step sizes for the explicit Euler method $\hat{E}$, the version $\hat{I}$ of the implicit Euler method and the Euler-exponential method $\hat{Z}$*

| $M$ | **2000** | **4000** | **8000** | **16000** |
|---|---|---|---|---|
| $\varepsilon(\hat{E}, M)$ | 46.6545 | 46.7107 | 46.6381 | 23.5562 |
| $\Delta(\hat{E}, M)/2$ | 0.023207 | 0.13302 | 0.31203 | 0.28929 |
| $\varepsilon(\hat{I}, M)$ | 6.6179 | 5.5739 | 3.9754 | 2.5181 |
| $\Delta(\hat{I}, M)/2$ | 0.030248 | 0.045375 | 0.037721 | 0.054798 |
| $\varepsilon(\hat{Z}, M)$ | 0.33533 | 0.2236 | 0.11426 | 0.037446 |
| $\Delta(\hat{Z}, M)/2$ | 0.066711 | 0.059289 | 0.077666 | 0.098786 |

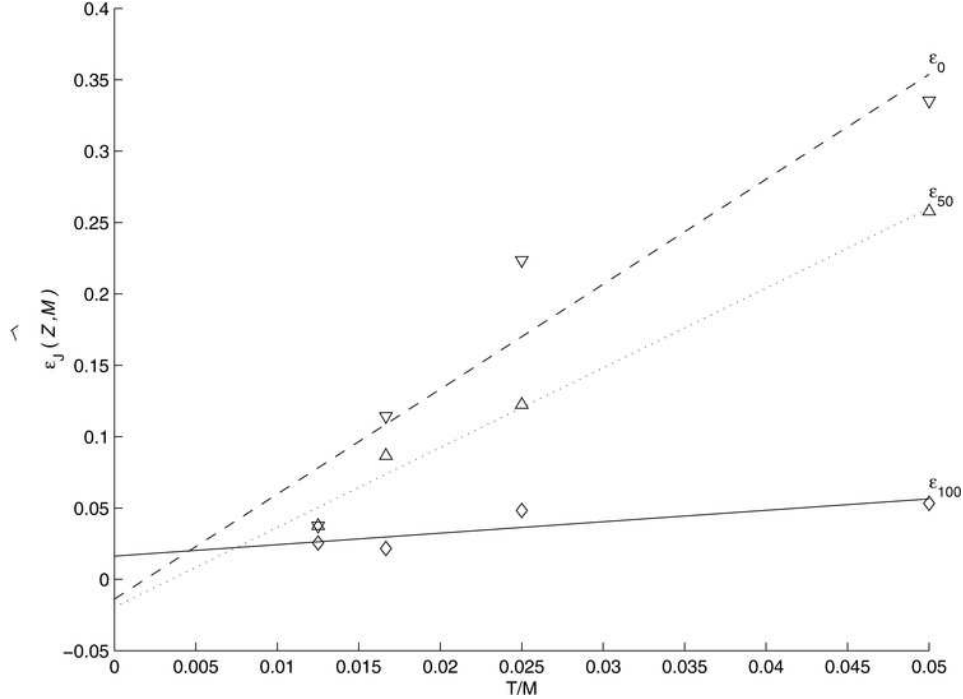

Fig. 2. *Errors versus step size for Scheme* 2: $J = 0$, $\triangledown$, $J = 50$, $\triangle$ *and* $J = 100$, $\diamond$.

It can be seen from both Figure 1 and Table 1 that the Euler-exponential scheme presents a superior performance than the other two numerical methods for this example. For instance, the error induced by Scheme 3 with $M = 16 \cdot 10^3$ is substantially greater than the error induced by Scheme 2 with $M = 2 \cdot 10^3$. Furthermore, the accuracy of Scheme 2 is very good for large time step sizes. This suggests that $\hat{Z}$ shows great promise for the long time integration of stochastic Schrödinger equations.

Finally, Figure 2 shows precision-step size diagrams. In particular, this figure gives the errors $\varepsilon_J(\hat{Z}, M)$, with $J = 0, 50, 100$, versus the step size $T/M$. Moreover, it presents the best least square linear approximation of each $\varepsilon_J(\hat{Z}, \cdot)$. From Figure 2, we see that the errors induced by $\hat{Z}$ closely follow a straight line. This is in a good agreement with Theorem 2.

**Acknowledgments.** The author is greatly indebted to Professor R. Rebolledo and Professor D. Talay for many helpful comments and suggestions. The author wishes to thank the Pontificia Universidad Católica de Chile, where part of this research was carried out, for financial support and hospitality.

DEPARTAMENTO DE INGENIERÍA MATEMÁTICA
UNIVERSIDAD DE CONCEPCIÓN
CASILLA 160 C, CONCEPCIÓN
CHILE
E-MAIL: cmora@ing-mat.udec.cl
URL: www.ing-mat.udec.cl/~cmora/